\documentclass{amsart}
\usepackage{amscd}

\usepackage{xypic}
\input xy
\xyoption{all}
\usepackage{graphicx}
\def\Dot{\lower.2pc\hbox to 2pt{\hss$\bullet$\hss}}
\def\Circ{\lower.2pc\hbox to 2pt{\hss$\circ$\hss}}

\newcommand{\bbC}{{\mathbb C}}
\newcommand{\bbP}{{\mathbb P}}
\newcommand{\bbZ}{{\mathbb Z}}
\newcommand{\bbN}{{\mathbb N}}

\newcommand{\Or}{{\mathcal O}}
\newcommand{\Si}{\Sigma}

\newtheorem{lemma}{Lemma}[section]
\newtheorem{prop}[lemma]{Proposition}
\newtheorem{thm}{Theorem}
\newtheorem{other}[lemma]{Theorem}

\newtheorem{cor}[lemma]{Corollary}
\newtheorem{defn}{Definition}

\begin{document}
\title{The Orevkov invariant of an affine plane curve}
\author{Walter D. Neumann}
\address{Department of Mathematics\\Columbia University\\
NY 10027\\USA}
\email{neumann@math.columbia.edu}
\author{Paul Norbury}
\address{Department of Pure Mathematics\\
Adelaide University\\Australia 5005}
\email{pnorbury@maths.adelaide.edu.au}

\begin{abstract}
We show that although the fundamental group of the complement of an
algebraic affine plane curve is not easy to compute, it possesses a
more accessible quotient, which we call the Orevkov invariant.
\end{abstract}

\maketitle

\section{Introduction.}
An interesting topological invariant of an algebraic affine plane
curve $\Si\subset\bbC^2$ is its link at infinity $L_{\Si}$
obtained by intersecting $\Si$ with a large sphere in $\bbC^2$.  It
is quite easily computed using \cite{NeuCom,NeuIrr}.

From $L_{\Si}$ and the number of reducible components of $\Si$ one can
retrieve the arithmetic genus of $\Si$ \cite{NeuCom,NeuIrr}, and much
about the topology of the defining polynomial $f$ given by
$\Si=\{f=0\}$ including the topology of its generic fibre and
information about its local and global monodromy representations
\cite{DNeMon,NNoVan}.  The link $L_{\Si}$ also enables one to obtain
geometric information about the curve beyond homological information
\cite{NNoUnf}.  This raises the question: {\em how much can $L_{\Si}$
tell us about $\pi_1(\bbC^2-\Si)$?}

If $\Si$ is a generic fibre of its defining polynomial $f$ then
$L_{\Si}$ gives the topology of $\Si$ as an embedded curve in $\bbC^2$
and hence it determines $\pi_1(\bbC^2-\Si)$ (in fact,
$\pi_1(\bbC^2-\Si)=\bbZ$ in this case by Oka \cite{OkaTwo}.  See
Section~\ref{sec:nodal}, Theorem~\ref{th:oka}.)  In contrast, the
cuspidal and nodal curves $x^2=y^3$ and $x^2=y^3+1$ both have the
trefoil at infinity whereas the fundamental groups of their
complements are respectively the braid group $\{ a,b|aba=bab\}$ and
$\bbZ$.  This leads us to consider a common quotient of the different
fundamental groups.

\begin{defn}
  If $p$ is a singular point of a plane curve $\Si$, let
  $G_p=\pi_1(B_p\cap(\bbC^2-\Si))$ for a sufficiently small ball $B_p$
  around $p$.  The \emph{Orevkov invariant} of $\Si$ is
  \[\Or(\Si)=G/N\] where $G=\pi_1(\bbC^2-\Si)$ and $N$ is the normal
  closure in $G$ of the union over all singular points $p$ of $\Si$ of
  the images of the commutator subgroups $[G_p,G_p]$.
\end{defn}
{\em Remark.}  The image of $G_p$ in $G$ is only well-defined
up to conjugation but since we take the normal closure, $N$ is 
well-defined.

\vspace{6pt}
When the invariant $\Or(\Si)$ is abelian it is
$H_1(\bbC^2-\Si)=\bbZ^k$, where $k$ is the number of irreducible
components of $\Si$.  This invariant was first studied by Orevkov in
\cite{OreCom} who proved that it is abelian when $L_{\Si}$ is a
positive braid.  (He also showed that the invariant is abelian on
projective plane curves which generalises the (solved) Zariski
conjecture \cite{OreFun}.)  The invariant ignores the complications of
configurations of hyperplanes since it is abelian in these cases.

In this paper we will show that the Orevkov invariant of an algebraic
affine plane curve is invariant under a special class of deformations
of the curve.  
By ``deformation'' we would normally mean a complex
analytic family, but the deformations we consider need only be complex
analytic near the singularities.
\begin{defn}
  A deformation of a singularity germ is \emph{local component
  preserving} (LCP) if it does not change the number of local irreducible
  components near the singularity. An LCP deformation of a curve is a
  deformation that is an LCP deformation in a neighbourhood of each
  singularity, and is a smooth proper isotopy outside smaller 
  neighbourhoods of the singularities.
\end{defn}

This definition applies to projective or affine curves. For affine
curves the condition that the isotopy is proper
means that the link at infinity is
preserved.
For example, if $\Si$ is the image of a non-constant analytic map
$f_0\colon\bbC\to\bbC^2$ then any complex analytic deformation $f_t$
of $f_0$ which preserves the link at infinity gives an LCP deformation
of $\Si$. Examples of non-LCP deformations are:
\begin{enumerate}
\item 
the deformation
$\Si_\epsilon=\{(x,y):xy=\epsilon\}$ of $\Si_0=\{(x,y):xy=0\}$ (not LCP
at the singularity);
\item
the deformation $\Si_\epsilon=\{(x,y):xy^2+y=\epsilon\}$ of
$\Si_0=\{(x,y):xy^2+y=0\}$ (not proper; link at infinity is not preserved).
\end{enumerate}

\begin{thm} \label{th:inv}
    The Orevkov invariant is invariant under LCP deformations.
\end{thm}

A {\em polynomially parametrised curve} is a curve with irreducible
components given by algebraic maps from $\bbC$ to $\bbC^2$.  Each
irreducible component can be parametrised by a pair of single variable
polynomials and is thus a rational curve with one place at infinity.

Polynomially parametrised curves arise naturally as the set of
non-generic regular values of a self-map of $\bbC^2$ as follows. 
Consider a polynomial map $F:\bbC^2\rightarrow\bbC^2$ with
two-dimensional image.  The pre-image of a generic value consists of
$d$ points, say.  There are two classes of non-generic values: those
that are regular and those that are irregular.  The pre-image of a
non-generic regular value consists of fewer than $d$ points.  The
closure (in $\bbC^2$) of the set of non-generic regular values is a
polynomially parametrised curve.  This is because the missing points
in the pre-image lie ``at infinity''.  If we compactify $\bbC^2$ by
including a divisor at infinity so that $\bar{F}$ is well-defined
there, then the points of the divisor at infinity that map to finite
values under $\bar{F}$ consist of some of the rational curves of the
divisor at infinity.  These rational curves are disjoint and each
ratonal curve consists of exactly one point that is mapped to
infinity, and hence the image of the union of the ratonal curves is a
polynomially parametrised curve.  The image of the rational curves at
infinity might intersect the set of irregular values and hence it
gives the closure (in $\bbC^2$) of the non-generic regular points.

If one restricts to deformations of a polynomially parametrised curve
that remain within the class of polynomially parametrised curves, then
automatically the LCP condition is satisfied at singularities.  Hence,
from Theorem~\ref{th:inv}, we see that in many cases the Orevkov
invariant of a polynomially parametrised curve is related to its link
at infinity.

\begin{cor}
    The Orevkov invariant of a polynomially parametrised curve depends
    only on its link at infinity when the moduli space of polynomially
    parametrised curves with that link at infinity is connected.
\end{cor}

Let $X_{\Si}$ be a compactification of $\bbC^2$ on which $\Si$ meets
the divisor at infinity $D$ transversally.  The link at infinity
$L_{\Si}$ encodes the minimal such divisor $D$ (see, e.g.,
\cite{NeuIrr}) and it also encodes the canonical class of $X_{\Si}$
supported on $D$.  Each component of the link corresponds to a
rational curve in $D$, known as a {\em horizontal curve}.  When the
canonical class at each link component is negative enough on each
horizontal curve we can prove that the moduli space of polynomially
parametrised curves with link at infinity $L_{\Si}$ is connected and
use this to prove the following theorem.

\begin{thm}   \label{th:suff}
    For sufficiently negative canonical class on each horizontal curve
    of $X_{\Si}$, the Orevkov invariant $\Or(\Si)$ of a polynomially
    parametrised curve $\Si$ is abelian.
\end{thm}

We state precisely how negative is sufficient in
Section~\ref{sec:solve}, Theorem~\ref{th:suffp}, and in
Section~\ref{sec:pos} we show that the class of positive braids
$L_{\Si}$ is strictly contained in this set.  Note that Orevkov
\cite{OreCom} proved that the Orevkov invariant is abelian on any
curve with positive braid at infinity, not just polynomially
parametrised curves with such a property.  Since our sufficiently
negative condition is only a mild improvement of Orevkov's result in
the polynomially parametrised case, it is the connectivity of certain
moduli spaces proved in this paper that is the more significant
improvement of his result.

An example of a link at
infinity that satisfies the conditions of Theorem~\ref{th:suff} is the
$(5,2)$ cabling on the $(2,3)$ torus knot, which is not a positive
braid, and represented by the following splice diagram: 
$$
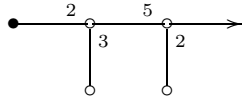
\objectmargin{0pt}\spreaddiagramrows{-3pt}
\spreaddiagramcolumns{3pt}\diagram
\Dot\rline^(.75){2}&\Circ\dline^(.25){3}\rline^(.75){5}&\Circ
\dline^(.25){2}\rto&\\ &\Circ&\Circ& \enddiagram
$$

\section{Nodal curves.}  \label{sec:nodal}
When the curve $\Si$ is a nodal curve---its only singularities are
simple double points---then ${\mathcal
O}(\Si)=\pi_1(\bbC^2-\Si)$.  The fundamental group of the
complement of a nodal plane curve has a long history beginning with
the Zariski conjecture, proven by Deligne \cite{DelGro} and Fulton
\cite{FulFun}, that the complement of a projective plane nodal curve
has abelian fundamental group.  This was generalised by Nori
\cite{NorZar} to projective and affine curves in smooth surfaces.

\begin{other}[Nori]    \label{th:nori}
Let $\Si$ and $D$ be curves on a smooth projective surface $X$, that
intersect transversally.  Assume that $\Si$ is nodal.  Denote the
number of singular points on a curve $C$ by $r(C)$.  Assume that
$C\cdot C>2r(C)$ for every irreducible component $C$ of $\Si$.  Then,
if $N$ denotes the kernel of $\pi_1(X-\Si\cup D)\rightarrow\pi_1(X-D)$,
$N$ is a finitely generated abelian group and the centraliser of $N$ 
is a subgroup of finite index.
\end{other}
In particular, when $\Si$ is a nodal curve in $\bbC^2$, and $X$ the
blow-up of $\bbP^2$ that resolves the singularities of $\Si$ at
infinity, and $D=X-\bbC^2$ the divisor at infinity, if
$\bar{\Si}_i\cdot\bar{\Si}_i>2r(\Si_i)$ for each irreducible component
$\Si_i$ of $\Si$ then $\pi_1(\bbC^2-\Si)$ is abelian (since
$X-D=\bbC^2$ and $N=\pi_1(\bbC^2-\Si)$ in the theorem.)

As mentioned in the introduction, when $\Si$ is the generic fibre
of its defining polynomial, $L_{\Si}$ determines
$\pi_1(\bbC^2-\Si)$.  
\begin{other}[Oka \cite{OkaTwo}]  \label{th:oka}
If $\Si$ is the generic fibre of a polynomial $f$ then
$\pi_1(\bbC^2-\Si)=\bbZ$.
\end{other}
\begin{proof}
    Choose a disk $D_s\subset\bbC$ to contain all of the atypical
    values of $f$ and a small disjoint disk $D_c$ that gives a
    neighbourhood of $c=f(\Si)$.  Join the disks by a path $\gamma$
    and put $D'=D_s\cup\gamma\cup D_c$.  Then each of the following
    arrows is a homotopy equivalence: \[ f^{-1}(D_s)\hookrightarrow
    f^{-1}(D_s\cup\gamma)\hookrightarrow
    f^{-1}(D')\hookrightarrow\bbC^2\] \[ f^{-1}(D'-\{
    c\})\hookrightarrow\bbC^2-\Si.\] Apply the Seifert Van-Kampen
    theorem to \[f^{-1}(D'-\{ c\})=f^{-1}(D_s\cup\gamma)\cup
    f^{-1}(D_c-\{ c\}).\]  Since $f^{-1}(D_c-\{ c\})\cong(D_c-\{
    c\})\times\Sigma\sim S^1\times\Si$ we get that $\pi_1(\bbC^2-\Si)$
    is the quotient of
    $\pi_1(\bbC^2)*\pi_1(S^1\times\Si)=\pi_1(S^1\times\Si)$ by the
    normal closure of $\pi_1(\Sigma)$, and thus
    $\pi_1(\bbC^2-\Si)=\bbZ$.
\end{proof}
{\em Remark}.  It is amusing to note that Theorem~\ref{th:oka} is
almost a consequence of Theorem~\ref{th:nori} which requires
$\bar{\Si}\cdot\bar{\Si}>0=2r(\Si)$ to conclude that
$\pi_1(\bbC^2-\Si)$ is abelian.  The self-intersection number of the
generic fibre of a polynomial is non-negative.  It is given as a sum
over linking numbers of components of $L_{\Si}$.  In \cite{NeuIrr} it
was shown that the polynomial is irregular at a point at infinity
precisely when the linking number for the corresponding link component
is $0$.  In particular, for a good polynomial, one which is never
irregular at infinity, the self-intersection number of the generic
fibre is strictly positive and Theorem~\ref{th:nori} implies
Theorem~\ref{th:oka}.  In fact, Theorem~\ref{th:nori} applies in most
cases since it is quite rare that the generic fibre of a polynomial
has zero self-intersection, or equivalently that it is irregular at
infinity at each link component of $L_{\Si}$.  Russell's bad field
generator \cite{RusGoo} and each polynomial in Kaliman's
classification of rational polynomials with a $\bbC^*$ fibre
\cite{KalRat} are examples of polynomials whose generic fibre has zero
self-intersection number, and hence Theorem~\ref{th:nori} does not
imply Theorem~\ref{th:oka} in these cases.

\vspace{6pt}
This paper will consist mainly of examples of curves with abelian
Orevkov invariant.  Two examples with {\em non-abelian} invariant are
as follows.

(i) For $f:\bbC^2\rightarrow\bbC$ the curve $\Si=f^{-1}(A)$ for a
finite set of at least two points $A\subset\bbC$ has non-abelian
Orevkov invariant.

(ii) The following example was constructed by Mutsuo Oka.  Let
\[ f(x,y)=(x^2-y)^2y-4x(x^2-y)+4=0.\] 
Then $\Si=f^{-1}(0)$ is a smooth irreducible curve with
\[\pi_1(\bbC^2-\Si)=\{ a,b|aba=bab\}=\Or(\Si).\]

\section{Invariance under local deformations.}  \label{sec:inv}
In this section we will prove Theorem \ref{th:inv} by showing that it
reduces easily to an argument of Orevkov's.

\begin{proof}
  Given a plane curve $\Si$, we can enclose its singularities in small
  balls $B_i$ and compute $\pi_1(\bbC^2-\Si)$ by the Van Kampen
  theorem in terms of the fundamental groups of the complement of
  $\Si$ in each of the balls as well as their exterior
  $X=\bbC^2-\bigcup_i B_i$. Under an LCP deformation, we can assume
  that the topology of $X-(\Si\cap X)$ does not change, so it suffices to show
  that the Orevkov invariant of each $(B_i, B_i\cap\Si)$ also does not
  change.
  
  We therefore restrict to the local situation, in which we have a
  germ of a singularity $(\Si,p)$ and a small ball $B$ around $p$ such
  that its boundary and the boundary of all smaller balls centered at
  $p$ intersect $\Si$ transversally. We take a LCP deformation that
  changes the intersection of $\Si$ with $\partial B$ only by an
  isotopy. In this situation we must show that the Orevkov invariant
  of $(B, B\cap\Si)$ does not change. We may change coordinates so
  that $p$ is the origin and $\Si$ is given by an equation $f(x,y)=0$
  of degree $n$ with $f(x,y)=x^n+a_{n-1}(y)x^{n-1}+\dots+a_0(y)$. We
  may then take the ball $B$ to be polydisk
  $\{(x,y):|x|\le\epsilon_1,|y|\le\epsilon_2\}=D_1\times D_2$, chosen
  such that $\Si$ intersects $\partial B$ only in the portion
  $\partial D_1\times D_2$.  Finally, we may assume our coordinates
  are chosen so that the set of $x$ for which $\Si\cap B$ contains
  more than one point $(x,y)$ whose $y$-coordinates have the same real
  part is in generic position in the sense of Orevkov's papers
  \cite{OreFun} and \cite{OreCom}, both before and after our
  deformation.  We denote this set $$P:=\{x\in D_1:\exists~ y_1\ne
  y_2\text{ with }(x,y_1)\in\Si, (x,y_2)\in\Si, \Re y_1=\Re y_2\}.$$
  Generic position implies that $P$ consists of a collection of paths
  oriented by increasing $|y_1-y_2|$.  Before deforming $\Si$, these
  paths start at the image of the one singular point and fan out to
  $\partial D_1$, since we choose the ball small enough that the only
  singular point of the projection $B\cap\Si\to D_1$ is the singular
  point of $\Si$.  After the deformation there may be several singular
  points, but we may assume that the behavior of the paths near
  $\partial D_1$ has changed only by an isotopy.  The paths then still
  start at singular values of the projection $B\cap \Si\to D_1$ and
  end at the boundary of the disk.  Moreover, away from singular
  points of the projection $B\cap\Si\to D_1$, the paths may now
  intersect, but Orevkov's genericity condition says they intersect at
  most in pairs or in triples, and these intersections are transverse
  (the triple intersections arise when three real parts of $y$-values
  coincide at a time).  The projection $B\cap\Si\to D_1$ is a covering
  map on the inverse image of $D_1-P$, of degree $d$, say, and, as in
  \cite{OreFun}, one can give a presentation of $\pi_1(B-\Si)$ with
  $d$ generators for each component of $D_1-P$, and with relations
  associated to the singular points of $B\cap \Si\to D_1$ and to the
  intersection points of paths comprising $P$.  The same inductive
  argument that Orevkov uses, moving along the paths starting from the
  singular points, then shows that the Orevkov invariant is free
  abelian with one generator for each irreducible component of $B\cap
  \Si$.
\end{proof}

\section{Splice diagrams.}
A convenient way to represent the link at infinity of an affine plane
algebraic curve uses the splice diagram of the link.  The splice
diagram of a link at infinity can be viewed as a cabled torus knot in
$S^3$, or the Puiseux expansion of the affine curve at infinity given
by expressing one coordinate in terms of the other, or the plumbing
graph of an efficient resolution at infinity of the affine curve.  We
will demonstrate these three views by focusing on the example of a
curve with knot at infinity.  This includes irreducible polynomial
curves.  The splice diagram of such a curve is given by
\begin{figure}[ht]
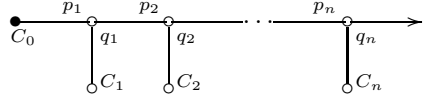

$$
\objectmargin{0pt}\spreaddiagramrows{-3pt}
\spreaddiagramcolumns{3pt}\diagram
\Dot\rline_(.1){C_0}^(.75){p_1}&\Circ\dline^(.25){q_1}^(.9){C_1}\rline^(.75)
{p_2}&\Circ\dline^(.25){q_2}^(.9){C_2}\rline&\dots\rline^(.75){p_n}&\Circ
\dline^(.25){q_n}^(.9){C_n}\rto&\\&\Circ&\Circ&&\Circ& 
\enddiagram
$$
\caption{Splice diagram at infinity.}
\label{fig:splice}
\end{figure}
where the pairs of integers $(p_k,q_k)$ are known as weights of the
splice diagram and the $C_k$ are simply labels for the valency one
nodes.  

The splice diagram represents a knot obtained by cabling a $(p_2,q_2)$
torus knot on a $(p_1,q_1)$ torus knot and then cabling a $(p_3,q_3)$
torus knot on that, and so on.  The earlier torus knots are
represented as ``virtual'' link components in the diagram by the
valency one nodes $C_k$.  For example, the components $C_0$ and $C_1$
form a Hopf link and the component $C_2$ is a $(p_1,q_1)$ torus knot
cabled on $C_1$.  The arrowhead represents the actual knot, a
$(p_n,q_n)$ cabling on the virtual knot $C_n$.

The linking number of any two (virtual) components can be calculated
from the splice diagram.  It is the product of all weights incident to
a path joining the two nodes.  Thus, $l(C_j,C_k)=p_jq_{j+1}\dots q_k$
and if we denote the knot (represented by the arrowhead) by $K$, then
$l(C_k,K)=p_kq_{k+1}\dots q_n$.

Not all knots arise as the knot at infinity of a plane curve.
Necessary and sufficent conditions on the weights of a splice diagram
of a knot at infinity are:
\begin{equation}  \label{eq:weight}
    \left\{\begin{array}{l}
    {\rm i}.\ p_k>1,\ q_k>1,\ (p_k,q_k)=1\\
    {\rm ii}.\ \Delta_k=p_k-p_{k-1}q_{k-1}q_k<0\quad\mbox{and}\quad
    \Delta_1=p_k-q_k<0\\
    {\rm iii}.\ p_{k+1}\in\bbN q_1q_2\dots q_k\oplus\bbN p_1q_2
    \dots q_k\oplus\dots\oplus\bbN p_{k-1}q_k\oplus\bbN p_k
    \end{array}\right.
\end{equation}
The third rather deep property involving linking numbers of virtual
link components is the semi-group condition of Abhyankar and Moh
\cite{AMoSem}.

A plane curve $\Si$ which has this knot at infinity is a degree
$q_1\dots q_n$ curve and furthermore we can choose coordinates $x$
and $y$ so that the defining polynomial $P(x,y)$ of $\Si$ has a single
highest degree term $y^{q_1q_2\dots q_n}$ and the highest degree
monomial in $x$ is $x^{p_1q_2\dots q_n}$.  Then one can expand $y$ in
terms of $x$:
\begin{equation}   \label{eq:puis} 
    y=x^{p_1/q_1}(a_{10}+a_{11}x^{e_{11}/q_1}+a_{12}x^{e_{12}/q_1}+...+
    x^{e_2/q_1q_2}(a_{20}+a_{21}x^{e_{21}/q_1q_2}+..
\end{equation}
for decreasing exponents
$0>e_{11}/q_1>e_{12}/q_1>...>e_2/q_1q_2>(e_2+e_{21})/q_1q_2>..$ where
$e_k=\Delta_k$ and the highest degree term in which the denominator of
the exponent does not divide $q_1\dots q_{k-1}$ is
$a_{k0}x^{e_k/q_1q_2...q_k}$, ($a_{k0}\neq 0$, $(e_k,q_k)=1$.)  We can
express $y$ more neatly as
\[ y=x^{p_1/q_1}(r_1(x^{-1/q_1})+x^{\Delta_2/q_1q_2}(r_2(x^{-1/q_1q_2})+..
+x^{\Delta_k/q_1..q_k}(r_k(x^{-1/q_1..q_k})+..\] 
for polynomials $r_k$ with $r_k(0)\neq 0$ and degree
$r_k<-\Delta_{k+1}/q_{k+1}$.  Alternatively, in keeping with the
approaches of \cite{ENeThr,NeuIrr} where the expansion at infinity is
obtained by compactifying and taking a local expansion around the
point at infinity, we can make the expansion homogeneous by
introducing the coordinate at infinity $z$.  The expansion becomes: 
\[ y=(x^{-1}z)^{p_1/q_1}z(r_1((x^{-1}z)^{1/q_1}) 
+..+(x^{-1}z)^{-\Delta_k/q_1..q_k}
(r_k((x^{-1}z)^{1/q_1..q_k})+..\]
and setting $x=1$, we get an expansion for $y$ in terms of $z$:
\[ y=z^{1-p_1/q_1}(r_1(z^{1/q_1})+z^{-\Delta_2/q_1q_2}
(r_2(z^{1/q_1q_2})+..+z^{-\Delta_k/q_1..q_k}(r_k(z^{1/q_1..q_k})+..\]
which gives the splice diagram
$$
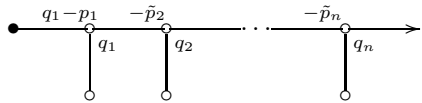
\objectmargin{0pt}\spreaddiagramrows{-3pt}
\spreaddiagramcolumns{3pt}\diagram
\Dot\rline^(.75){q_1-p_1}&\Circ\dline^(.25){q_1}\rline^(.75){-\tilde{p}_2}&
\Circ\dline^(.25){q_2}\rline&\dots\rline^(.75){-\tilde{p}_n}&\Circ\dline^(.25)
{q_n}\rto&\\&\Circ&\Circ&&\Circ& 
\enddiagram
$$
for $\tilde{p}_k=p_k+\xi_k$ where $\xi_k$ is chosen so that the edge
determinant is $-\Delta_k$.  By \cite{NeuCom} this local splice
diagram converts to the splice diagram at infinity given in (A).

The splice diagram also gives an efficient resolution around the point
of $\Si$ at infinity.  Compactify $\bbC^2$ to $\bbP^2$ and let
${\mathcal P}(y,z)=0$ define $\overline{\Si}$ in a neighbourhood of
the point $[1:0:0]\in\bbP^2$ (which we have supposed to be the point
where $\overline{\Si}$ meets the line at infinity.)  Each
characteristic pair in the splice diagram encodes a multiple blow-up
of the resolution of the point at infinity.  To resolve at the $k$th
step, one replaces $(y,z)$ by
$(y^az^{-\Delta_k}+Q_k(y^bz^{q_k}),y^bz^{q_k})$ where $a,b>0$ are
chosen so that $aq_k+b\Delta_k=1$.  The single variable polynomial
$Q_k$ is uniquely determined.  The curve $\overline{\Si}$ is
``resolved'' at infinity when it meets the divisor at infinity
transversally.  Each node $C_k$ represents a curve in the divisor at
infinity and the defining polynomial of the curve $\overline{\Si}$ in
a neighbourhood of $C_k$ is the polynomial obtained at the $k$th step
of the resolution.

The approximate roots of the defining polynomial $P(x,y)$ of $\Si$
\[ P_2=P^{1/q_2\dots q_n},\dots,
P_k=P^{1/q_k\dots q_n},\dots,P_n=P^{1/q_n}\] 
are each defined uniquely by the respective condition
\[ \deg_y(P-P_k^{q_k\dots q_n})<q_1\dots q_n-q_1\dots q_{k-1}.\]
It is easy to calculate each of these by setting 
\[ P_k=y^{q_1\dots q_{k-1}}+b_1(x)y^{q_1\dots q_{k-1}-1}+\dots+
b_{q_1\dots q_{k-1}}(x)\] and solving for $b_i(x)$.  (We define
$P_1=P^{1/q_1\dots q_n}$ by $\deg_y(P-P_1^{q_1\dots q_n})<q_1\dots
q_n-1$ and {\em a priori} $P_1=y-b(x)$, although by choice of
coordinates $P_1=y-y_0$ for a constant $y_0$.)

The zero set of the polynomial $P_k$ defines a curve that meets the
splice diagram, or equivalently the divisor at infinity, at $C_k$. 
This follows from the fact that when we resolve $P$, to get
$\tilde{P}$, the approximate $1/q_k\dots q_n$ root $P_k$ resolves to
give the approximate $1/q_k\dots q_n$ root of $\tilde{P}$. 
Furthermore, restricted to $C_k$, $\tilde{P}$ is a degree $q_k\dots
q_n$ polynomial with approximate $1/q_k\dots q_n$ root given by the
resolution of $P_k$ so the latter must be a coordinate.

Consider all holomorphic functions on $\Si$, meromorphic at infinity. 
The orders of the poles of these functions form a semi-group in $\bbN$
known as the Weierstrass semi-group, with complement the Weierstrass
gap sequence.  Using Riemann-Roch, the size of the Weierstrass gap
sequence can be shown to be equal to the arithmetic genus $g$ of $\Si$
which is easily calculated in terms of the splice diagram:
\begin{equation}  \label{eq:genus}
1-2g=q_1..q_n+\sum_{i=1}^{n-1}p_i(1-q_i)q_{i+1}..q_n+p_n(1-q_n).
\end{equation}
\begin{other}  \label{th:wegen}
The Weierstrass semi-group is generated by the orders of the poles of
the holomorhic functions on $\Si$ given by the coordinates $x$ and
$y$ and the approximate roots $P_2,\dots,P_n$.
\end{other}
\begin{proof}
The semi-group, $H_n$, generated by the poles of $x$ and $y$ and the
approximate roots $P_2,\dots,P_n$ at $\infty\in\overline{\Si}$ is
contained in the Weierstrass semi-group so it is sufficient to show
that the size of the gap $G_n=\bbN-H_n$ is equal to the size of the
Weierstrass gap sequence.

The order of the pole of the approximate root $P_k$ can be calculated
from the link at infinity.  It is given by the linking number of the
knot at infinity of $\Si$ with the knot at infinity of $P_k$ and the
latter is the virtual knot represented by the node $C_k$ of the splice
diagram.  Thus the order of the pole of $P_k$ is
$l(C_k,K)=p_kq_{k+1}\dots q_n$, and $l(C_n,K)=p_n$,
$l(C_0,K)=q_1q_2\dots q_n$, $l(C_1,K)=p_1q_2\dots q_n$, are the
respective orders of the poles of $P_n$, $x$ and $y$.  For
$k=1,\dots,n$ put \[ H_k=\bbN q_1q_2\dots q_k\oplus\bbN p_1q_2\dots
q_k\oplus \bbN p_2q_3\dots q_k\oplus\dots\oplus\bbN
p_{k-1}q_k\oplus\bbN p_k\] so the semi-group condition of Abhyankar
and Moh can be restated as $p_k\in H_{k-1}$.

Any element $a\in H_k$ has a unique normal form given by
\[a=\sum_{i=0}^na_k\cdot l(C_i,C_k)=a_0\cdot q_1q_2\dots 
q_k+\sum_{i=1}^{k-1}a_i\cdot p_iq_{i+1}\dots q_k+a_k\cdot p_k\] 
for non-negative integers $a_i$ satisfying $a_i<q_i$ ($i=1,..,k$)
and $a_0$ unrestricted.  This follows easily from the fact that since
$p_kq_k\in H_{k-1}q_k\subset H_k$, in any expression for $a$ we can
reduce the coefficient of $p_k$ by multiples of $q_k$, redistributing
the quantity amongst the previous coefficients of $l(C_i,C_k)$ for
$i<k$.  A similar step can then be taken to adjust the coefficient of
$p_{k-1}q_k$ by multiples of $q_{k-1}$ to ensure it is less than
$q_{k-1}$, and so on.

There are two further steps to the proof.  We must first find numbers
$m_k$ for $k=1,...,n$ such that $a>m_k\Rightarrow a\in H_k$, and then
use the unique normal form to count $I_k=\#H_k\cap [0,m_k]$ and hence
calculate the size of the gap $G_k=\bbN-H_k$.  Note that the number
$m_k$ need not be optimal.
\begin{lemma}  \label{th:mk}
$m_k=q_km_{k-1}+p_k(q_k-1)$
\end{lemma}
\begin{lemma}  \label{th:ik}
$I_k=q_kI_{k-1}+\frac{1}{2}(q_k-1)(p_k-1).$
\end{lemma}
We will delay the proofs of these two results.  Thus, 
\[|G_k|=m_k+1-I_k=q_k|G_{k-1}|+\frac{1}{2}(q_k-1)(p_k-1)\]
and the size of the Weierstrass gap, which is the arithmetic genus
given in (\ref{eq:genus}), satisfies the same recursion relation.
Thus it is sufficient to prove the theorem for a single Puiseux pair
$(p_1,q_1)$.

Since $(p_1,q_1)=1$, given $1\leq r\leq q_1$ we can choose $1\leq
a<q_1$ and $1\leq b<p_1$ such that $ap_1-bq_1=r$.  Then
\begin{eqnarray*}
    p_1q_1+r&=&p_1q_1+ap_1-bq_1\\
    &=&(p_1-b)q_1+ap_1\in H_1
\end{eqnarray*}
for each $r=1,..,q_1$ and hence for any $r\geq 1$, so $m_1=p_1q_1$.
(In fact, the same argument shows we can use $m_1=p_1q_1-p_1-q_1$, but
that is not necessary here.)

And $I_1=\#\{a_0q_1+a_1p_1\leq p_1q_1\}=(1/2)(p_1+1)(q_1+1)$ by
counting integer lattice points under $q_1x+p_1y=p_1q_1$.  (The value
$p_1q_1$ is represented twice, and only counted once.)

Hence $|G_1|=p_1q_1+1-(1/2)(p_1+1)(q_1+1)=(1/2)(p_1-1)(q_1-1)$ which
is the arithmetic genus and hence equals the size of the Weierstrass
gap.
\end{proof}

\begin{proof} of Lemma~\ref{th:mk}.
Since $q_kH_{k-1}\subset H_k$, then $q_k(m_{k-1}+1)\in H_k$.  Any
$a>p_kq_k-p_k-q_k$ lies in the semi-group generated by $p_k$ and
$q_k$, hence any $a>q_k(m_{k-1}+1)+p_kq_k-p_k-q_k$ lies in $H_k$. 
Thus $m_k=q_km_{k-1}+p_k(q_k-1)$.
\end{proof}

\begin{proof} of Lemma~\ref{th:ik}
Since $m_k=q_km_{k-1}+p_k(q_k-1)$, and since we may choose $a_k<q_k$,
we get
\begin{eqnarray*}
I_k&=&\#\left\{\sum_{i=0}^ka_i\cdot p_iq_{i+1}..q_k\leq m_k\right\}\\
&=&\#\left\{\sum_{i=0}^{k-1}a_i\cdot p_iq_{i+1}..q_{k-1}\cdot q_k\leq 
q_km_{k-1}\right\}\ \ \ ({\rm when}\ a_k=q_k-1)\\
&+&\#\left\{\sum_{i=0}^{k-1}a_i\cdot p_iq_{i+1}..q_{k-1}\cdot q_k\leq 
q_km_{k-1}+p_k\right\}\ \ \ (a_k=q_k-2)\\
&\vdots&\\
&+&\#\left\{\sum_{i=0}^{k-1}a_i\cdot p_iq_{i+1}..q_{k-1}\cdot q_k\leq 
q_km_{k-1}+(q_k-1)p_k\right\}\ \ \ (a_k=0)\\
&=&q_kI_{k-1}+\left[\frac{p_k}{q_k}\right]+\left[\frac{2p_k}{q_k}\right]
+...+\left[\frac{(q_k-1)p_k}{q_k}\right]\\
&=&q_kI_{k-1}+\frac{1}{2}(q_k-1)(p_k-1).
\end{eqnarray*}
The $i$th of the $q_k$ terms contributes $I_{k-1}$ plus the number of
multiples of $q_k$ less than $(i-1)p_k$, given by $[ip_k/q_k]$, the
greatest integer part of $ip_k/q_k$, and hence the second last
expression $q_kI_{k-1}+...$ follows.  The sum of greatest integer
parts is obtained by counting integer lattice points under the graph
$p_kx=q_ky$.
\end{proof}

\section{Polynomially parametrised curves.}
If $\Si$ is a polynomially parametrised curve then any deformation
within the class of polynomially parametrised curves preserves the
irreducible components globally and hence locally.  Hence any such
deformation that also preserves the link at infinity is an LCP
deformation and we can apply Theorem~\ref{th:inv} to show that the
Orevkov invariant of $\Si$ depends only upon its connected component 
of polynomially parametrised curves with a given link at infinity.

In this section we will describe the moduli space of polynomially
para\-metrised curves.  This will allow us to study the connected
components of the moduli space and in particular when such moduli
spaces are connected.

\subsection{Moduli spaces.}

We will begin with a description of the moduli space of {\em
irreducible} polynomially para\-metrised curves with a given splice
diagram at infinity.

\begin{other}  \label{th:moduli}
Let the pair of polynomials $(x(t),y(t))$ define a reduced rational
curve $\Si$.  Then the following two properties are equivalent:

(A) the link at infinity of $\Si$ has splice diagram 
$$
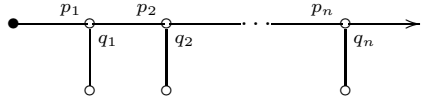
\objectmargin{0pt}\spreaddiagramrows{-3pt}
\spreaddiagramcolumns{3pt}\diagram
\Dot\rline^(.75){p_1}&\Circ\dline^(.25){q_1}\rline^(.75){p_2}&\Circ
\dline^(.25){q_2}\rline&\dots\rline^(.75){p_n}&\Circ\dline^(.25){q_n}\rto&\\
&\Circ&\Circ&&\Circ& 
\enddiagram
$$

(B) The defining polynomial for $\Si$ is given by
$P(x,y)=y_{n+1}(x,y_1,y_2,..,y_n)$ where $y_1=y$ and for each
$k=1,...,n$, 
\begin{equation}   \label{eq:yk} 
y_{k+1}(x,y_1,y_2,..,y_k)=y_k^{q_k}+\sum_e
a_ex^{e_0}y_1^{e_1}y_2^{e_2}..y_k^{e_k}
\end{equation}
where $e=(e_0,..,e_k)$ is summed over all tuples such that the
order of the pole at $\infty$ (the $t$-degree) of the corresponding
monomial is no greater than the order of the pole at $\infty$ of
$y_k^{q_k}$.  As a polynomial in $t$,
$\deg_ty_k=p_kq_{k+1}q_{k+2}\dots q_n$ for $(p_k,q_k)=1$ and
$y_{n+1}\equiv 0$.
\end{other}

\begin{proof}
$(A)\Rightarrow(B)$ 

Given a curve parametrised by $(x(t),y(t))$, there is an algorithm to
determine the polynomial which vanishes on the curve, $P(x,y)=0$.  The
splice diagram at infinity guides the algorithm and the holomorphic
functions $\{y_2,\dots,y_n\}$ with the stated properties are produced
along the way.  The algorithm runs as follows.

Take a linear combination $y_{10}=y^{q_1}+\lambda x^{p_1}$ so that
$deg_ty_{10}<p_1q_1$.  By iteratively adding monomials
$y_{1,k+1}=y_{1k}+\alpha_{ij}x^iy^j$ we can get rid of successively
lower powers of $t$, reducing the degree $deg_ty_{1,k+1}<deg_ty_{1k}$.
This procedure stops when we can no longer reduce the degree using
monomials in $x$ and $y$ and produces $y_2(x,y)$ with pole at infinity
of order $\deg_ty_2=p_2q_3\dots q_n$ for the following reason.

The order of the pole of $y_2(x,y)$ lies in the Weierstrass semi-group
$W$ of the curve $\Sigma$.  By Theorem~\ref{th:wegen}, this is the
semi-group generated by the orders of the poles of $x$ and $y$ and the
approximate roots.  Use the ordering of the approximate roots to
define sub-semi-groups 
\begin{equation}   \label{eq:subsemi}
    W_1\subset W_2\subset\ldots W_n=W
\end{equation}
where $W_k$ is given by those elements of $W$ generated by
$q_1q_2\dots q_n$, $p_1q_2\dots q_n$, $p_2q_3\dots q_n$, \ldots,
$p_kq_{k+1}\ldots q_n$.  In terms of $H_k$ defined in the proof of
Theorem~\ref{th:wegen}, $W_k=q_nq_{n-1}\ldots q_{k+1}H_k$.  The
polynomial $y_2(x,y)$ behaves in a similar way to the approximate root
$P_2$, in that it meets the splice diagram at infinity below the
second Puiseux pair.  This is because it meets the splice diagram at
infinity before the third Puiseux pair since it possesses only a
single Puiseux pair itself, so the order of its pole lies in $W_2$,
and the order of its pole does not lie in $W_1$ since then we could
reduce the degree of $y_2(x,y)$ further by a polynomial in $x$ and $y$
and we should not have stopped the procedure.

By the semi-group condition of Abhyankar and Moh, $p_2$ lies in the
semi-group generated by $p_1$ and $q_1$, or equivalently one can
reduce the $t$-degree of $y_2^{q_2}$ by subtracting a monomial in $x$
and $y$.  As before, we iteratively reduce the $t$-degree of
$y_2^{q_2}$ using monomials in $x$, $y$ and $y_2$, and produce $y_3$
when the procedure stops.  The same argument as above shows that
$deg_ty_3=p_3q_4..q_n$ and the algorithm continues.  Finally we are
left with $y_{n+1}$ which consists of powers of $t$ that lie in the
Weierstrass gap.  But then $y_{n+1}\equiv 0$ since no holomorphic
functions on $\Si$ can have such poles.

Note that the $y_k$'s are {\em not} the approximate roots of the
defining polynomial of the curve.  In fact, the polynomial $y_k$ is
not well-defined since there are relations amongst the variables
$x,y_1,y_2,...$ and also since $y_k$ can be adjusted by low degree
terms.\\
\\ 
$(B)\Rightarrow(A)$ The expansion of $y$ in terms of $x$ is most
easily seen by putting $x=w^{-q_1\dots q_n}$ so $w$ is a branch of
$x^{-1/q_1\dots q_n}$.  This allows us to solve for
\[ t=a_{-1}w^{-1}+a_0+a_1w+a_2w^2+\dots,\ a_{-1}\neq 0\]
and then
\[ y=w^{-p_1q_2\dots q_n}(r_1(w^{q_2\dots q_n}) 
+\dots +w^{-e_kq_{k+1}\dots q_n}(r_k(w^{q_{k+1}\dots q_n})+\dots\]
where $w^{-e_kq_{k+1}\dots q_n}$ is the first term with exponent not
divisible by $q_kq_{k+1}\dots q_n$ and $q_{k+1}\dots q_n$ is defined
by requiring that $(e_k,q_k)=1$.  Note that an expansion for $y$ in
terms of $w$ always exists in such a form but {\em a priori}
$e_kq_{k+1}\dots q_n$ does not necessarily equal $\Delta_kq_{k+1}\dots
q_n$.

Now suppose there exists $y_2(x,y)$ satisfying B1, B2 and B3.  Then we
can express $y_2=w^{-p_2q_3\dots q_n}\eta(w)$ for $\eta$ holomorphic
and $\eta(0)\neq 0$.  The leading (most negative) power of $w$ is the
same as the leading power of $w$ not divisible by $q_2\dots q_n$.
Hence it comes from the term $y^{q_1}$.  Thus
\[ -p_2q_3\dots q_n=-q_1p_1q_2q_3\dots q_n-e_2q_3\dots q_n.\]
Any other term $x^iy^j$ does not contribute because it gives a more
positive power since at best one can get
\[ iq_1\dots q_n+jp_1\dots q_n-e_2q_3\dots q_n<q_1p_1\dots q_n
-e_2q_3\dots q_n\] since $iq_1+jp_1<p_1q_1$.

More generally, $y_k=w^{-p_kq_{k+1}\dots q_n}\eta(w)$ for $\eta$
holomorphic and $\eta(0)\neq 0$.  The leading power of $w$ not
divisible by $q_kq_{k+1}\dots q_n$ comes from the term
$y_{k-1}^{q_{k-1}}$ so
\begin{equation}   \label{eq:degw} 
-p_kq_{k+1}\dots q_n=q_{k-1}p_{k-1}q_k\dots q_n-e_kq_{k+1}\dots q_n
\end{equation}
and no other term $y_{k-1}^{i_{k-1}}\dots y_1^{i_1}x^{i_0}$
contributes since it can at best give the weighted degree of the
monomial minus $e_kq_k\dots q_n$ and the weighted degree of the
monomial is strictly less than $q_{k-1}p_{k-1}q_k\dots q_n$.

From (\ref{eq:degw}) we see that
$e_k=\Delta_k=p_k-p_{k-1}q_{k-1}q_k<0$ by (\ref{eq:weight}) ii.
\end{proof}
Here is an example to demonstrate the algorithm.\\
\\
{\em Example.} 
Let $x(t)=t^{12}+t$ and $y(t)=t^8+t^2$.  The algorithm of 
the previous proof enables us to produce the defining polynomial of 
$x$ and $y$ and to prove that its splice diagram at infinity is
$$
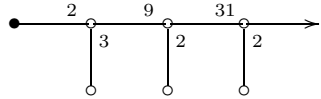
\objectmargin{0pt}\spreaddiagramrows{-3pt}
\spreaddiagramcolumns{3pt}\diagram
\Dot\rline^(.75){2}&\Circ\dline^(.25){3}\rline^(.75){9}&\Circ
\dline^(.25){2}\rline^(.75){31}&\Circ\dline^(.25){2}\rto&\\
&\Circ&\Circ&\Circ&
\enddiagram
$$
\begin{itemize}
    \item[Step 1.] $deg_t\{x^2-y^3\}=18$ so set $y_2(t)=x^2-y^3$;
    \item[Step 2.] $deg_t\{(x^2-y^3)^2-9x^3\}=31$ so set 
                   $y_3(t)=(x^2-y^3)^2-9x^3$;
    \item[Step 3.] $deg_t\{((x^2-y^3)^2-9x^3)^2+(16/3)y\cdot y_2^3\}=61$; 
    \item[Step 4.] add to 3. a multiple of $x\cdot y_2\cdot y_3$ to get 
                   rid of the $t^{61}$ term;
    \item[Step 5.] continue this to get rid of $t^{60},t^{59},\dots,t^0$
                   and find that only powers of $t$ that lie in the 
		   semi-group generated by 12, 8, 18 and 31 arise;
    \item[Step 6.] thus, $P(x,y)=((x^2-y^3)^2-9x^3)^2+(16/3)y(x^2-y^3)^3+...$
\end{itemize}
Alternatively, in step 2, $deg_t\{(x^2-y^3)^2-9xy^3\}=31$ so we may
set $y_3(t)=(x^2-y^3)^2-9xy^3$ and continue with the algorithm.
Although the process is not unique, and hence $P$ is not unique as a
polynomial in $x,y,y_2$ and $y_3$, as a polynomial in $x,y$, $P(x,y)$ is
unique.\\

The example demonstrates both the working of the algorithm and the
fact that the $y_j$ are not unique, in contrast with the unique
approximate roots.\\

The conditions (B) in Theorem~\ref{th:moduli} lead to an explicit
expression for the moduli space of polynomially parametrised curves
with given splice diagram at infinity as an algebraic variety.  For
$x(t)=\sum a_it^i$ and $y=\sum b_jt^j$ the variety is given as a set
of polynomial equations in the coefficients $\{ a_i,b_j\}$ as follows.

The expression for $y_{k+1}$ in (\ref{eq:yk}) {\em a priori} has
$t$-degree $p_kq_kq_{k+1}..q_n$ but since we require that $y_{k+1}$
has $t$-degree $p_{k+1}q_{k+2}q_{k+3}..q_n<p_kq_kq_{k+1}..q_n$, by
(\ref{eq:weight}) ii, then each of the coefficients of $t^m$ for
$p_kq_kq_{k+1}..q_n\geq m>p_{k+1}q_{k+2}q_{k+3}..q_n$ must vanish. 
For $m\in W_k$, defined in (\ref{eq:subsemi}), the coefficient of
$t^m$ is canceled by a monomial in $x$ and $y_i$, $i\leq k$.  For each
$m\notin W_k$ , define $T_{m,k}(\{ a_i,b_j\})$ to be the coefficient
of $t^m$ in $y_k$.  We call the variety $\{T_{m,k}=0\}$ the moduli
space of polynomially parametrised curves with given splice diagram at
infinity.

Each polynomial $T_{m,k}$ has the property that it only depends on
$b_j$ and $a_i$ for $i\geq e(k)-m$, where $e(1)=p_1q_1q_2\ldots q_n$
and $e(k)=e(k-1)+(q_k-1)p_kq_{k+1}\ldots q_n$, and if equality occurs
then that coefficient appears in $T_{m,k}$ linearly with non-zero
coefficient.  By non-zero coefficient, we mean that it is a 
polynomial in the $\{ b_j\}$ that is required to be non-zero (earlier 
in the algorithm.)

When the link at infinity has more than one component, or equivalently
the polynomially parametrised curves are reducible, if the splice
components meet only on the root vertex, or equivalently the curve
components do not intersect at infinity, then the moduli space is
simply a product of the moduli spaces for each component, minus a
divisor along which the curves meet at infinity.  In general, any two
irreducible components of the link at infinity agree along an initial
iterated cabling and this introduces a further equation giving
equality of the initial coefficients of the respective pairs of
polynomials.

It is often difficult to tell when the variety is connected or even
when it is non-empty.  In the case of a knot at infinity, the
semi-group theorem of Abhyankar and Moh \cite{AMoSem} and its converse
by Sathaye and Stenerson \cite{SStPla} give necessary and sufficient
conditions for the moduli space of curves (not necessarily rational)
to be non-empty.

\subsection{Solving the equations.}   \label{sec:solve}
In the previous section we constructed a variety that gives the moduli
space of an irreducible component of a polynomially parametrised curve
with a given link at infinity.  In this section we will be more
explicit, and solve the equations for a given link at infinity.  It is
important to note that we only consider splice diagrams that are links
at infinity of algebraic curves so $\Delta_k=p_k-p_{k-1}q_{k-1}q_k<0$
for all $k$.

We will begin with a short discussion of the canonical class divisor
associated to $L_{\Si}$.  The splice diagram for $L_{\Si}$ gives an
efficient resolution at infinity for $\Sigma$.  It gives the valency
$>2$ curves in the plumbing diagram of the divisor at infinity
$D\subset X_{\Si}$.  It can be arranged that the canonical class $K$
of $X_{\Si}$ is supported on the divisor $D$, so $K$ is given by a
sequence of multiplicities, one for each irreducible curve in $D$.  At
each component of $L_{\Si}$, the multiplicity of $K$ is given by the
local intersection $K\cdot\Sigma$ at that component.  Instead of
working with $K\cdot\Sigma$, we prefer to work with $-K\cdot\Sigma-1$
which we denote by $d$ at each component of the link, and more
generally at each virtual component of the link, or node of the splice
diagram.

\begin{figure}[ht]
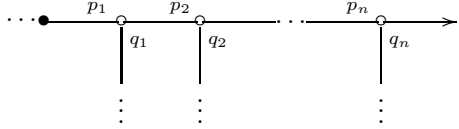

$$
\objectmargin{0pt}\spreaddiagramrows{-3pt}
\spreaddiagramcolumns{3pt}\diagram
\dots\Dot\rline^(.75){p_1}&\Circ\dline^(.25){q_1}\rline^(.75){p_2}&\Circ
\dline^(.25){q_2}\rline&\dots\rline^(.75){p_n}&\Circ\dline^(.25){q_n}\rto&\\
&\vdots&\vdots&&\vdots& 
\enddiagram
$$
\caption{Path in splice diagram at infinity.}
\label{fig:path}
\end{figure}

The number $d$ can be calculated at a node in the splice diagram of
$L_{\Si}$ by taking a path from the node to the root---the node that
represents the proper transform of $\bbP^1=\bbP^2-\bbC^2$---of the
splice diagram.  Such a path is represented by Figure~\ref{fig:path}
where other arrows may exist in the diagram at the places marked with
dots.  If we notate $d$ at the $k$th node by $d_k$ then $d_k$ is
calculated recursively by
\begin{equation}
d_k=q_kd_{k-1}+\Delta_k,\ \  d_1=p_1+q_1.
\end{equation}

The following theorem is a more precise statement of the sufficiently
negative condition of Theorem~\ref{th:suff}.
\begin{other}  \label{th:suffp}
When $d_n\geq p_1q_2\dots q_n -1$ for each component of $L$ then the
moduli space of polynomially parametrised curves $\Si$ with
$L_{\Si}=L$ is connected.
\end{other}
\begin{proof}
To start, consider each irreducible component of $\Si$ separately and
hence use the description of the moduli space given in the previous
section.  A straight-forward way to solve the equations is to
arbitrarily choose the coefficients of $y(t)$, and solve for the
coefficients of $x(t)$, and use this to show that when the moduli
space is non-empty it is connected.  As described in the properties of
the $T_{m,k}$, the equations $T_{m,k}=0$ form a diagonal system in which each
coefficient of $x$ appears for the first time linearly with non-zero
coefficient.

For each $k=2,\dots,n$, the equation $\deg_ty_k=p_kq_{k+1}q_{k+2}\dots
q_n$ gives rise to 
\begin{equation}  \label{eq:count}
p_{k-1}q_{k-1}q_kq_{k+1}q_{k+2}\dots
q_n-1-p_kq_{k+1}q_{k+2}\dots q_n=-\Delta_kq_{k+1}q_{k+2}\dots q_n-1
\end{equation}
vanishing coefficients since {\em a priori} the degree of $y_k$ is
$p_{k-1}q_{k-1}q_kq_{k+1}q_{k+2}\dots q_n-1$.  The coefficients of
powers of $t$ that do not lie in the sub-semi-group $W_k$ give rise to
the equations $T_{m,k}=0$ and these are easily solved, as a diagonal
system when the number of coefficients of $x$ is no less than the
number of equations.  (In fact, many coefficients of $x$ are wasted 
since we need only solve $T_{m,k}=0$ for $m\notin W_k$.)  We need:  
\begin{equation}  \label{eq:suffp}
    q_1q_2\dots q_n\geq\sum_{k=2}^n\Delta_kq_{k+1}q_{k+2}\dots q_n+1.
\end{equation}
Put \[a(n)=(q_1+p_1)q_2\dots q_n+\sum_{k=2}^n\Delta_kq_{k+1}q_{k+2}\dots 
q_n-d_n\]
so (\ref{eq:suffp}) becomes 
\[ d_n\geq p_1q_2\dots q_n-1-a(n).\]
Using $\Delta_n-d_n=-d_{n-1}q_n$ we see that
$a(n)=a(n-1)q_n$ or $a(k)=a(k-1)q_k$ for all $k$ and since $a(2)=0$ by
induction $a(n)=0$.  Hence, $d_n\geq p_1q_2\dots q_n-1$.

The proof is not complete since the vanishing of the specified
coefficients only guarantees $\deg_ty_k\leq p_kq_{k+1}\dots q_n$ and
not $\deg_ty_k=p_kq_{k+1}\dots q_n$.  If, by good fortune, we do get
equality each time, then we can calculate the dimension of the space
of solutions to the system as the degree of $y(t)$ plus the number of
coefficients of $x(t)$ unused when solving the system.
\begin{lemma}  \label{th:reg}
The moduli space of polynomially parametrised curves with given link
at infinity is either of the expected dimension and connected, or
empty.
\end{lemma}
\begin{proof}
We have arranged that the coefficient of $t^{p_kq_{k+1}q_{k+2}\dots
q_n+1}$ vanishes.  The coefficient of $t^{p_kq_{k+1}q_{k+2}\dots q_n}$
is a polynomial $B_k$ in the coefficients of $y(t)$.  One of two cases
must occur.  Either:

(i) the polynomial $B_k$ is identically zero thus the moduli space of
polynomially parametrised curves with that particular link at infinity
is empty; or

(ii) the polynomial $B_k$ is not identically zero thus the
coefficients of $y(t)$ that give a smaller degree for $y_k$ form a set
of codimension $\geq 1$ and the moduli space of coefficients of $y(t)$
that give the right degree for $y_k$ is a non-empty connected space.
\end{proof}

If each component of $L_{\Si}$ gives rise to a non-empty moduli space
of polynomially parametrised curves, then by construction we can fit
the moduli spaces of the components of $\Si$ together to get the full
moduli space.  Two link components agree along an initial path of
cablings, or equivalently the respective compactification divisors at
infinity agree on an initial path of blow-ups, when the same top
coefficients of $y(t)$ are chosen along with the same polynomials
$y_k$.  Again, accidental agreement along a further path of cablings
occurs either on a codimension $\geq 1$ set or the moduli space of
polynomially parametrised curves with given link at infinity is empty.

This completes the proof of the theorem. 
\end{proof}
{\em Remarks} (i) We can ensure that each non-empty moduli space
contains an immersion.  When an unused coefficient $a_i$ of $x$
appears for the first time, and hence linearly, in $T_{m,k}$, say, we
can use $a_i$ as a parameter.  The polynomial $x(t)$ then depends on
the parameter $a_i$ and $(x'(t),y'(t))\neq (0,0)$ for generic choice
of $a_i$.

(ii) In the next section we prove that for positive braids
Lemma~\ref{th:reg} can be strengthened to show that the space of
solutions is always non-empty.

We are now almost in a position to prove Theorem~\ref{th:suff}.  The
following lemma proves that a particular class of immersions has
abelian fundamental group and hence abelian Orevkov invariant.

\begin{lemma}  \label{th:imm}
Let $\Si=\sqcup S_i\rightarrow\bbC^2$ be an immersion.  If
$2g(\Si_i)-2>K.\bar{\Si_i}$ on each component $\Si_i$ then
$\pi_1(\bbC^2-\Si)$ is abelian.
\end{lemma}
\begin{proof}
This is a simple application of Nori's theorem.  First, assume that the
immersion of $\Si$ is nodal.  Then on each component
\[\bar{\Si_i}\cdot\bar{\Si_i}+2-2p_a(\Si_i)=-K\cdot\bar{\Si_i}\]
and since the arithmetic genus is given by
$p_a(\Si_i)=g(\Si_i)+r(\Si_i)$ where $r(\Si_i)$ is the number of
double points of $\Si_i$ then we have
\[\bar{\Si_i}\cdot\bar{\Si_i}-2r(\Si_i)=-K.\bar{\Si_i}+2g-2\]
and when this is strictly positive, for each component $\Si_i$, Nori's
theorem applies to conclude that $\pi_1(\bbC^2-\Si)$ is abelian.

When the immersion of $\Si$ is not nodal, blow up triple points and
higher multiplicity points.  Add to the divisor at infinity the
exceptional curves whose union is given by $E$ so that Nori's theorem
is applied to $\tilde{Si}\rightarrow X-D\cup E$.  Then after each
blow-up, the number $\bar{\Si_i}\cdot\bar{\Si_i}$ (use $\Si$ to denote
its proper transform also) is reduced by the same as the number
$-K.\bar{\Si_i}$, which maintains the inequality, and the arithmetic
genus is possibly reduced, which improves the inequality.  Thus, the
general case follows from the nodal case.
\end{proof}

\begin{proof} of Theorem~\ref{th:suff}.

By the invariance of $\Or(\Si)$ under LCP deformations and the
connectedness of the moduli space of polynomially parametrised curves
under LCP deformations, it is sufficient to calculate $\Or(\Si)$ for
one curve $\Si$ with given link at infinity.

As explained in Remark (ii) after Theorem~\ref{th:suffp}, each moduli
space contains an immersion.  When each component is rational with one
point at infinity, the condition $2g(\Si_i)-2>K.\bar{\Si_i}$ in
Lemma~\ref{th:imm} becomes $-K.\bar{\Si_i}-1>1$ and the left hand side
is the number $d$ we associate to the link component of the splice
diagram.  The condition $d>1$ is certainly fulfilled on each component
since $d_n\geq p_1q_2\dots q_n -1>1$.  Hence $\Or(\Si)$ is abelian
and this is true on the entire moduli space, so the Theorem is proven.
\end{proof}

\subsection{Positive braids.}   \label{sec:pos}
The link at infinity of an affine algebraic curve inherits a natural
braid structure.  In fact there is a braid structure for each line in
$\bbC^2$ obtained by projection along the line.  The condition that
the splice diagram of a link at infinity gives a positive braid is a
condition on the weights $p_i,q_i$ of each branch of the splice
diagram given in Figure~\ref{fig:path}.  

In order to make the largest possible set of links at infinity
positive braids, we choose to project along the $x$ direction (assume
degree $x(t)>$ degree $y(t)$).  In the following proposition, if we
were to use the braid structure from the projection along the $y$
direction, then the condition $d_n>p_1q_2\dots q_n$ would be replaced
by $d_n>q_1q_2\dots q_n>p_1q_2\dots q_n$ which is a stronger condition
and hence fewer links would be positive braids.

\begin{prop}  \label{th:posb}
    The link at infinity of an affine curve is a positive braid
    precisely when $d_n>p_1q_2\dots q_n$ along each path of the splice
    diagram shown in Figure~\ref{fig:path}.
\end{prop}
\begin{proof}
Each component of a positive braid is again a positive braid so we
will restrict to each path in the splice diagram.  We will actually
prove that $d_k>p_1q_2\dots q_k$, $k=1,\dots,n$.  Thus we claim that
the corresponding virtual link components are also positive braids.

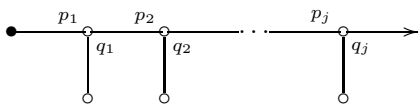
\begin{figure}[ht]
$$
\objectmargin{0pt}\spreaddiagramrows{-3pt}
\spreaddiagramcolumns{3pt}\diagram
\Dot\rline^(.75){p_1}&\Circ\dline^(.25){q_1}\rline^(.75){p_2}&\Circ\dline^(.25)
{q_2}\rline&\dots\rline^(.75){p_j}&\Circ\dline^(.25){q_j}\rto&\\
&\Circ&\Circ&&\Circ& 
\enddiagram
$$
\label{fig:splice3}
\end{figure}
\begin{lemma}
Using the framing supplied by the braid projection, the self-linking
number of the braid in Figure~\ref{fig:splice3} is
\[ l_j=p_jq_j+p_1q_2\dots q_j-d_j.\] 
Equivalently, it is the braid index of the braid.
\end{lemma}
\begin{proof}
Prove this by induction.  To begin, the $(p_1,q_1)$ torus knot has
self-linking number $l_1=p_1q_1-q_1=p_1q_1+p_1-d_1$, as required.  To
see the recursive relationship for $l_j$, suppose we cable $q_{j+1}$
parallel strings along Figure~\ref{fig:splice3} using the framing
supplied by the braid projection.  Then the self-linking number, or
braid index, of this new braid is given by $q_{j+1}^2l_j$.  Now
consider a $q_{j+1}$ cabling on Figure~\ref{fig:splice3} with weight
$p_{j+1}$.  Note that $p_{j+1}$ gives the linking number of (the
virtual copy of) Figure~\ref{fig:splice3} with each of the $q_{j+1}$
cables and a cabling already links $q_{j+1}l_j$ times with (the
virtual copy of) Figure~\ref{fig:splice3}.  Thus, the new self-linking
number is 
\begin{equation}   \label{eq:rec}
l_{j+1}=q_{j+1}^2l_j+(q_{j+1}-1)(p_{j+1}-q_{j+1}l_j)=q_{j+1}l_j+
q_{j+1}p_{j+1}-p_{j+1}.
\end{equation} 
In particular, $l_j-p_jq_j$ and $p_1q_2\dots q_j-d_j$ satisfy the same
recursion relation.  Since $l_1-p_1q_1=p_1-d_1$, the result follows.
\end{proof}
Now, a positive braid arises when the $p_{j+1}-q_{j+1}l_j$ extra
twists in (\ref{eq:rec}) is positive.  Hence $p_k>l_{k-1}q_k$,
$k=1,\dots,n$ and $d_k>p_1q_2\dots q_k$ for $k=1,\dots,n$.

Note, too, that the recursive relationship for $d_j$ shows that
$d_n>p_1q_2\dots q_n$ implies $d_k>p_1q_2\dots q_k$ for $k=1,\dots,n$.
\end{proof}
Since the condition in Proposition~\ref{th:posb} is stronger than that
in Theorem~\ref{th:suffp} we get the folllowing corollary.
\begin{cor}
The moduli space of polynomially parametrised curves with positive
braid at infinity is connected.
\end{cor}

\begin{prop}  \label{th:nonem}
The moduli space of polynomially parametrised curves with positive
braid at infinity is non-empty.
\end{prop}
\begin{proof}
In the proof of Theorem~\ref{th:suffp}, we produced points in the
variety for a given splice diagram by solving for the coefficients of
$f$ in terms of the coefficients in $g$.  The variety can be empty
when there is accidental vanishing of the coefficients of any of
$t^{p_kq_{k+1}..q_n}$ for $k=1,..,n$ (where $n$ is the number of
Puiseux pairs of the particular component.)  Since at each component
of a positive braid there are $n-1$ extra degrees of freedom,
\[d_n\geq p_1q_2\dots q_n + 1\geq p_1q_2\dots q_n -1\]
then these can be used to set the coefficients of $t^{p_kq_{k+1}..q_n}$
to be $1$.  Thus, for each component the moduli space is non-empty.

When we put the components together, components branching off at
different points in the splice diagram correspond to different choices
of coefficients that lie in the semi-group of the curve, and hence
waste no degrees of freedom.  Thus the moduli space of polynomial
curves with positive braid at infinity is non-empty.
\end{proof}

For positive braids, we can prove Theorem~\ref{th:suff} without using
Nori's theorem.  We used Nori's theorem to calculate the Orevkov
invariant for a special representative in the moduli space of
polynomially parametrised curves.  For positive braids, we can do such
a calculation explicitly.

\begin{other}   \label{th:pert}
    Given a curve $(F(t),G(t))$ define 
    \[(f(t),g(t))=(F(t^d)+\epsilon h(t),G(t^d)).\] 
    Then for most degree $h<d$ and for small enough $\epsilon$ the
    Orevkov invariant of $(f(t),g(t))$ is a quotient of the Orevkov
    invariant of $(F(t),G(t))$.
\end{other}
\begin{proof}
Choose a large ball $B_R$ that realises the link at infinity of
$(F(t),G(t))$.  Choose a path from $0$ to $\infty$ in $\bbC^*$ that
avoids self-intersections of $(F,G)$.  Choose $\epsilon$ small enough
so that a tube of radius $\epsilon R$ intersects the curve
$(F(t^d),G(t^d))$ only in an embedded strip in the curve.  Thus along
the path the $d$ deformations remain disjoint from the rest of
$(F,G)$.  The $d$ deformations may intersect each other.  The branch
point at $t=0$ gives the relation that each of the $d$ elements in the
fundamental group are equal.  Thus, the Orevkov invariant of the
deformed curve $(f(t),g(t))=(F(t^d)+\epsilon t,G(t^d))$ is a quotient
of the Orevkov invariant of $(F(t),G(t))$.  The Orevkov invariants are
the same if $B_R$ realises the link at infinity of $(f(t),g(t))$.
\end{proof}

Next we give a construction of a rational curve with positive braid at
infinity in order to apply Theorem~\ref{th:pert}.
\begin{other}  \label{th:pospert}
If $(F(t),G(t))$ defines the component of a positive braid 
$$
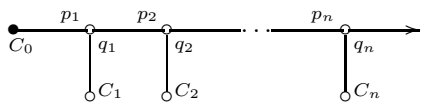
\objectmargin{0pt}\spreaddiagramrows{-3pt}
\spreaddiagramcolumns{3pt}\diagram
\Dot\rline_(.1){C_0}^(.75){p_1}&\Circ\dline^(.25){q_1}^(.9){C_1}\rline^(.75)
{p_2}&\Circ\dline^(.25){q_2}^(.9){C_2}\rline&\dots\rline^(.75){p_n}&\Circ
\dline^(.25){q_n}^(.9){C_n}\rto&\\&\Circ&\Circ&&\Circ& 
\enddiagram
$$ 
then for degree $h(t)<(d_n-q_1q_2\ldots q_n)q_{n+1}$, the curve
$$(f(t),g(t))=(F(t^{q_{n+1}})+\epsilon h(t),G(t^{q_{n+1}}))$$ 
defines the component
$$
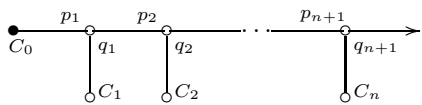
\objectmargin{0pt}\spreaddiagramrows{-3pt}
\spreaddiagramcolumns{3pt}\diagram
\Dot\rline_(.1){C_0}^(.75){p_1}&\Circ\dline^(.25){q_1}^(.9){C_1}\rline^(.75)
{p_2}&\Circ\dline^(.25){q_2}^(.9){C_2}\rline&\dots\rline^(.75){p_{n+1}}&\Circ
\dline^(.25){q_{n+1}}^(.9){C_n}\rto&\\&\Circ&\Circ&&\Circ& 
\enddiagram
$$ 
where any $p_{n+1}$ that satisfies
$p_{n+1}-p_nq_nq_{n+1}+d_nq_{n+1}=d_{n+1}>p_1q_2\dots q_nq_{n+1}$ 
is uniquely determined by the choice of $h$.
\end{other}
\begin{proof}
    When we solve the diagonal system of equations $T_{m,k}=0$ for a
    positive braid at infinity to get $(F(t),G(t))$, the coefficients
    of $F(t)$ appear in order from the highest power of $t$ to the
    lowest.  Furthermore, the last $d_n-q_1q_2\ldots q_n$ coefficients
    of $F(t)$ are irrelevant to the system of equations.  If the
    degree of $\epsilon h(t)$ is less than $(d_n-q_1q_2\ldots
    q_n)q_{n+1}$ then $\epsilon h(t)$ does not affect the equations
    $T_{m,k}=0$ for $k<n$, and features only in the equations
    $T_{m,n}=0$.  The equations $T_{m,n}=0$ form a diagonal sustem in 
    the coefficients of $\epsilon h(t)$ and hence we can solve for 
    any $p_{n+1}$ that gives rise to a positive braid, or equivalently 
    satsifies the inequality
    $p_{n+1}-p_nq_nq_{n+1}+d_nq_{n+1}=d_{n+1}>p_1q_2\dots q_nq_{n+1}$.
\end{proof}

Thus, by induction the Orevkov invariant of a polynomially
parametrised curve $\Si$ with positive braid at infinity is abelian
since it is a quotient of the Orevkov invariant of a simpler curve
obtained by reducing by one the number of Puiseux pairs on a component
of $\Si$.  Eventually the curve is reduced to a configuration of
lines, where the invariant is abelian.

\section{Further examples.}
 
The invariance of the Orevkov invariant under LCP deformations reduces
the problem of calculating the invariant, or searching for a curve
with non-abelian invariant, to understanding the finitely many
topological types at infinity of curves with given degree.  One can
enumerate polynomially parametrised curves by their splice diagrams at
infinity and the connected components of the moduli space of
polynomially parametrised curves with that splice diagram at infinity. 
When the curves are irreducible the list for degree up to $12$ is as
follows:

\begin{figure}[ht]
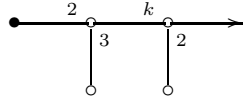

$$
\objectmargin{0pt}\spreaddiagramrows{-3pt}
\spreaddiagramcolumns{3pt}\diagram
\Dot\rline^(.75){2}&\Circ\dline^(.25){3}\rline^(.75){k}&\Circ
\dline^(.25){2}\rto&\\ &\Circ&\Circ& \enddiagram
$$
\caption{degree 6, $k=3,5,..,11$}
\end{figure}
\begin{figure}[ht]
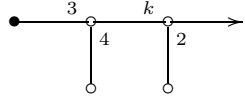

$$
\objectmargin{0pt}\spreaddiagramrows{-3pt}
\spreaddiagramcolumns{3pt}\diagram
\Dot\rline^(.75){3}&\Circ\dline^(.25){4}\rline^(.75){k}&\Circ
\dline^(.25){2}\rto&\\ &\Circ&\Circ& \enddiagram
$$
\caption{degree 8, $k=3,7,9,..,23$}
\end{figure}
\begin{figure}[ht]
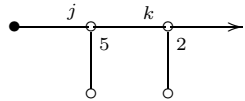

$$
\objectmargin{0pt}\spreaddiagramrows{-3pt}
\spreaddiagramcolumns{3pt}\diagram
\Dot\rline^(.75){j}&\Circ\dline^(.25){5}\rline^(.75){k}&\Circ
\dline^(.25){2}\rto&\\ &\Circ&\Circ& \enddiagram
$$
\caption{degree 10, $j=2,3,4$, odd $k=jm+5n<10j$}
\end{figure}
\begin{figure}[ht]
$$
\objectmargin{0pt}\spreaddiagramrows{-3pt}
\spreaddiagramcolumns{3pt}\diagram
\Dot\rline^(.75){2}&\Circ\dline^(.25){3}\rline^(.75){j}&
\Circ\dline^(.25){2}\rline^(.75){k}&\Circ\dline^(.25){2}\rto&
&&\Dot\rline^(.75){2}&\Circ\dline^(.25){3}\rline^(.75){l}&\Circ
\dline^(.25){4}\rto&\\
&\Circ&\Circ&\Circ& 
&&&\Circ&\Circ&\enddiagram
$$
\caption{degree 12, $j=3,5,..,11$, $k=j,..,4j-1$, 
$l=3,5,..,23$}
\label{fig:deg12}
$$
\objectmargin{0pt}\spreaddiagramrows{-3pt}
\spreaddiagramcolumns{3pt}\diagram
\Dot\rline^(.75){3}&\Circ\dline^(.25){4}\rline^(.75){j}&\Circ
\dline^(.25){3}\rto&
&&\Dot\rline^(.75){5}&\Circ\dline^(.25){6}\rline^(.75){k}&\Circ
\dline^(.25){2}\rto&\\
&\Circ&\Circ& 
&&&\Circ&\Circ&\enddiagram
$$
\caption{degree 12, $j=3,7,9,..,35$, $k=5,11,15,17,21,23,..,59$}
\end{figure}

In each of the above diagrams, by $a,..,b$ we mean all of the odd
integers from $a$ to $b$, except in the case of $k$ in
Figure~\ref{fig:deg12}, which misses some odd values.

The degree $<6$, 7 and 11 curves have positive braids at infinity and
we have not listed them nor positive braids for degree 6, 8, 10 and
12.  They are covered by Theorem~\ref{th:suff} along with many on the
list above (for example, when the degree is 6 only $k=3$ is not
covered.)  Using MAPLE for those curves above of degree $<12$ not
covered by Theorem~\ref{th:suff} we have shown that each moduli space
of polynomially parametrised curves is connected and that the Orevkov
invariant of each is abelian.

So far we have only been able to calculate the Orevkov invariant of a
curve by calculating the fundamental group of the complement of some
curve obtained by deforming the original curve and applying
Theorem~\ref{th:inv}.  Given the many cases where the Orevkov
invariant depends only on the link at infinity we would hope to be
able to calculate the invariant without calculating the fundamental
group of the complement of a curve.

One can construct smooth examples of curves $\bbC\rightarrow\bbC^2$
with non-abelian Orevkov invariant.  An example is as follows. 
Consider the representation of 
\[\pi_1(\bbC-{\rm 2\ points})=\{\gamma_1,\gamma_2\}\] 
into the braid group
\[B_3=\{\sigma_1,\sigma_2|\sigma_1\sigma_2\sigma_1=\sigma_2\sigma_1\sigma_2\}\]
given by 
\[\gamma_1\mapsto\sigma_2^{-2}\sigma_1\sigma_2^2\quad\mbox{and}\quad
\gamma_2\mapsto\sigma_1^2\sigma_2\sigma_1^{-2}.\]  
This gives rise to a smoothly embedded surface in $\bbC^2$ described
via two maps from a disk $\bbC$ to $\bbC$.  The first map, $x$, is a
3-fold branched cover with two branch points.  Take a regular point
$c$ for $x$ and a loop $\gamma$ starting and ending at $c$ that moves
around one of the branch points of $x$.  The pre-image
$x^{-1}(\gamma)$ is a path containg the three preimages $x^{-1}(c)$. 
The second map, $y$, is a 2-fold branched cover that wraps
$x^{-1}(\gamma)$ around its single branch point in such a way that it
has a double point.

The map $(x,y):\bbC\rightarrow\bbC^2$ defines a smoothly embedded
surface $\Si$ with \[\pi_1(\bbC^2-\Si)=\{ a,b|aba=bab\}.\] This cannot
be an algebraic curve since by the Abhyankar-Moh-Suzuki theorem any
embedded disk is equivalent to the standard disk.  We can add double
points to this curve, whilst preserving the Orevkov invariant, to make
it seem more like a polynomially parametrised curve.  Nevertheless,
the link at infinity of this does not arise as the link at infinity of
an algebraic affine plane curve.  It would be interesting to
characterise those links at infinity that do arise from a similar
construction of a non-abelian Orevkov invariant.

{\em Acknowledgements.}  The authors would like to acknowledge useful
discussions with S.Yu. Orevkov and Gavin Brown.


\begin{thebibliography}{99}

\bibitem{AMoSem} S. Abhyankar and T.T. Moh 
\emph{On the semigroup of a meromorphic curve},
Proc. Int Symp. Algebraic Geometry, Kyoto (1977), 249-414.

\bibitem{BroTop} S.A. Broughton
\emph{On the topology of polynomial hypersurfaces}
Proc. AMS Symp. Pure Math. {\bf 40} (1983), 167-178.

\bibitem{DelGro} Pierre Deligne 
\emph{Le groupe du complement d'une courbe plane n'ayant que des
points ordinaires est abelien (d'apres W. Fulton)},
Seminaire Bourbaki, Lect. Notes in Math. {\bf 842}, Springer Verlag (1981), 
1-10.

\bibitem{DNeMon} A. Dimca and A. Nemethi
\emph{On the monodromy of complex polynomials},
Duke Math. J. {\bf 108} (2001), 199-209.

\bibitem{ENeThr} Eisenbud, D. and Neumann, W.D.
\emph{Three-dimensional link theory and invariants of plane curve 
singularities.}
Ann. Math. Stud. {\bf 110}, Princeton Univ. Press (1985).

\bibitem{FulFun} W. Fulton
\emph{On the fundamental group of the complement of a node curve},
Ann. of Math. {\bf 111} (1980), 407-409.

\bibitem{KalRat} S. Kaliman
\emph{Rational polynomials with a $\bbC^*$-fiber},
Pacific J. Math. {\bf 174} (1996), 141-194.

\bibitem{LibAle} A. Libgober 
\emph{Alexander polynomial of plane algebraic curves and cyclic multiple 
planes}, Duke Math. J. {\bf 49} (1982), 833-851.

\bibitem{NeuCom} W.D. Neumann
\emph{Complex algebraic curves via their links at infinity},
Invent. Math. {\bf 3} (1989), 445-489.

\bibitem{NeuIrr} W.D. Neumann
\emph{Irregular links at infinity of complex affine plane curves}
Quarterly J. Math. {\bf 50} (1999), 301-320.

\bibitem{NNoVan} W.D. Neumann and P. Norbury
\emph{Vanishing cycles and monodromy of complex polynomials},
Duke Math. J. {\bf 101} (2000), 487-497.

\bibitem{NNoUnf} W.D. Neumann and P. Norbury
\emph{Unfolding polynomial maps at infinity},
Math. Ann. {\bf 318} (2000), 149-180.

\bibitem{NRuUnf} W.D. Neumann and L. Rudolph 
\emph{Unfoldings in knot theory}, 
Math. Ann. {\bf 278} (1987), 409-439 and Corrigendum {\bf 282} (1988),
349-351.

\bibitem{NorZar} Madhav V. Nori 
\emph{Zariski's conjecture and related problems},
Ann. Sci. Ecole Norm. Sup. (4) {\bf 16} (1983), 305-344.

\bibitem{OkaTwo} Mutsuo Oka 
\emph{Two transforms of plane curves and their fundamental groups},
J. Math. Sci. Univ. Tokyo {\bf 3} (1996), 399-433.

\bibitem{OreFun} S.Yu. Orevkov
\emph{The fundamental group of the complement of a plane algebraic curve},
Mat. Sb. {\bf 137 (179)} (1988), 267-277.

\bibitem{OreCom} S.Yu. Orevkov
\emph{The commutant of the fundamental group of the complement of a
plane algebraic curve}, 
Russian Math. surveys {\bf 45} (1990), 221-222.

\bibitem{RusGoo} P. Russell
\emph{Good and bad field generators},
J. Math. Kyoto Univ. {\bf 17} (1977), 319-331.

\bibitem{SStPla} A. Sathaye and J. Stenerson
\emph{On plane polynomial curves}
Algebraic geometry and its applications, C.L. Bajaj, Ed., 
Springer (1994), 121-142.

\bibitem{SerAlg} J.P. Serre
\emph{Algebraic groups and class fields},
Grad. Texts in Math. {\bf 117}, Springer-Verlag (1988).


\end{thebibliography}
\end{document}